\documentclass[11pt]{article}
\usepackage{geometry}
\usepackage{mathrsfs}
\usepackage{color}
\usepackage{epsfig}
\usepackage{graphicx}
\usepackage{epstopdf}
\usepackage{cite}
\usepackage{mathrsfs}
\usepackage{amsfonts}
\usepackage{amsmath}
\usepackage{amsfonts,amssymb,color}
\usepackage{dsfont}
\usepackage{curves}
\usepackage{mathrsfs}
\usepackage{pifont}
\usepackage{amssymb}
\usepackage{latexsym,amsmath,amssymb,amsfonts,epsfig,graphicx,cite,psfrag}
\usepackage{eepic,color,colordvi,amscd}
\def\qed{\hfill \rule{4pt}{7pt}}
\def\pf{\noindent {\it Proof. }}

\renewcommand{\baselinestretch}{1.2}

\title{A tight linear bound to the chromatic number of $(P_5, K_1+(K_1\cup K_3))$-free graphs\footnote{Supported by NSFC No. 11931106 and 12101117, and by NSFJS No. BK20200344}}
\author{Wei Dong$^{1,}$\footnote{Email: weidong@njxzc.edu.cn }, Baogang Xu$^{2,}$\footnote{Email: baogxu@njnu.edu.cn, or baogxu@hotmail.com.}, Yian Xu$^{3,}$\footnote{Email: yian$\_$xu@seu.edu.cn}
\\ $^{1}$\small School of Information and Engineering
\\ \small Nanjing Xiaozhuang University, Nanjing,
211171, China
\\ $^{2}$\small Institute of Mathematics, School of Mathematical Sciences
\\ \small Nanjing Normal University, Nanjing, 210023, China
\\$^{3}$\small School of Mathematics, Southeast University, 2 SEU Road, Nanjing, 211189, China}
\date{}
\newtheorem{theorem}{Theorem}[section]

\newtheorem{lemma}{Lemma}[section]

\newtheorem{prob}{Problem}[section]
\newtheorem{conj}{Conjecture}[section]

\def\qed{\hfill \rule{4pt}{7pt}}
\renewcommand{\baselinestretch}{1.15}
\begin{document}
\maketitle
\begin{abstract}
Let $F_1$ and $F_2$ be two disjoint graphs. The union $F_1\cup F_2$ is a graph with vertex set $V(F_1)\cup V(F_2)$ and edge set $E(F_1)\cup E(F_2)$,  and the join $F_1+F_2$ is a graph with vertex set $V(F_1)\cup V(F_2)$ and edge set $E(F_1)\cup E(F_2)\cup \{xy\;|\; x\in V(F_1)\mbox{ and } y\in V(F_2)\}$.  In this paper, we present a characterization to $(P_5, K_1\cup K_3)$-free graphs, prove that $\chi(G)\le 2\omega(G)-1$ if $G$ is $(P_5, K_1\cup K_3)$-free. Based on this result, we further prove that  $\chi(G)\le $max$\{2\omega(G),15\}$ if $G$ is a $(P_5,K_1+(K_1\cup K_3))$-free graph. We also construct a $(P_5, K_1+( K_1\cup K_3))$-free graph $G$ with $\chi(G)=2\omega(G)$.
\begin{flushleft}
\end{flushleft}
\end{abstract}

\section{Introduction}
All graphs considered in this paper are finite and simple.
Let $G$ be a graph. The vertex set of a complete subgraph of $G$ is called a {\em clique} of $G$, and the {\em clique number $\omega(G)$} of $G$ is the maximum size of cliques of $G$. We use $P_k$ and $C_k$ to denote a \textit{path} and a \textit{cycle} on $k$ vertices respectively.

Let $G$ and $H$ be two vertex disjoint graphs. The \textit{union} $G\cup H$ is the graph with $V(G\cup H)= V(G)\cup V(H)$
and $E(G \cup H) = E(G)\cup E(H)$. Similarly, the \textit{join} $G + H$ is the graph with $V(G + H) = V(G)\cup V(H)$ and $E(G + H) =
E(G) \cup E(H)\cup \{xy| \mbox{for each pair} \ x \in V(G) \mbox{ and } y \in V(H)\}$.

For a subset $X\subseteq V(G)$, let $G[X]$ denote the subgraph of $G$ induced by $X$. A \textit{hole} of $G$ is an induced cycle of length at least 4,  and a $k$-{\em hole} is a hole of length $k$. A $k$-{\em hole} is said to be an {\em odd $($even$)$ hole} if $k$ is odd (even). An \textit{antihole} is the complement of some hole. An {\em odd} (resp. {\em even}) antihole is defined analogously.

We say that $G$ induces $H$ if $G$ has an induced subgraph isomorphic to $H$,  and say that $G$ is $H$-{\em free} if $G$ does not induce $H$. Let $\mathcal{H}$ be a family of graphs. We say that $G$ is ${\cal H}$-free if $G$ induces no member of ${\cal H}$.

A coloring of $G$ is an assignment of colors to the vertices of $G$ such that no two adjacent
vertices receive the same color. The minimum number of colors required to color $G$ is called the {\em chromatic number} of $G$, and is denoted by $\chi(G)$. Obviously we have that $\chi(G)\ge \omega(G)$. However, determining the upper bound of the chromatic number of some family of graphs $G$, especially, giving a function of $\omega(G)$ to bound $\chi(G)$ is generally very difficult.  A family $\mathcal{G}$ of graphs is said to be $\chi$-{\em bounded} if there is a function $f$ such that
$\chi(G) \le f(\omega(G))$ for every $G\in \mathcal{G}$, and if such a function $f$ does exist to $\mathcal{G}$, then $f$ is said to be a {\em binding function} of $\mathcal{G}$ \cite{gyarfas1}. A graph $G$ is said to be \textit{perfect} if  $\chi(H)=\omega(H)$ for each induced subgraph $H$. Thus the binding function for perfect graphs is $f(x)=x$. The famous {\em Strong Perfect Graph Theorem} \cite{CRST06} states that a graph is perfect if and only if it is (odd hole, odd antihole)-free.  Erd\H{o}s \cite{E59} showed that for any positive integers $k$ and $\l$, there exists a graph $G$ with $\chi(G)\ge k$ and without cycles of length less than $\l$. This result motivates the study of the chromatic number of $\mathcal{H}$-free graphs for some $\mathcal{H}$. Gy\'{a}rf\'{a}s \cite{gyarfas1, G87}, and Sumner \cite{sumner1} independently, proposed the following conjecture.

\begin{conj}{\em \cite{G87, sumner1}}\label{Gyarfas87}
For every tree $T$, $T$-free graphs are $\chi$-bounded.
\end{conj}

Interested readers are referred to \cite{RS04,SS20,SR19} for more information on Conjecture~\ref{Gyarfas87} and related problems.
Gy\'{a}rf\'{a}s \cite{G87} proved that $\chi(G) \leq {(k-1)}^{\omega(G) - 1}$ for $k\ge 4$ if $G$ is $P_k$-free and $\omega(G)\ge 2$.  Then the  upper bound  was improved to ${(k-2)}^{\omega(G) - 1}$ by  Gravier {\em et al.} \cite{GHM2003}. The problem of determining whether
the class of $P_t$-free graphs ($t\ge5$) admits a polynomial $\chi$-binding function  remains open.

\begin{prob}{\em \cite{S16}}\label{Schiermeyer16}
Are there polynomial functions $f_{P_k}$ for $k\ge 5$ such that $\chi(G)\le f_{P_k}(\omega(G))$ for  every $P_k$-free graph $G$$?$
\end{prob}

Since $P_4$-free graphs are perfect, finding an optimal binding function for $P_5$-free graphs attracts much attention. Esperet {\it et al}. \cite{ELMM13} proved that $\chi(G)\leq 5\times 3^{\omega(G)-3}$ for $P_5$-free graphs.

\begin{theorem}\label{thm-esperet}{\em (\cite{ELMM13})}
$\chi(G)\le 5\cdot 3^{\omega(G)-3}$ for $P_5$-free graphs $G$ with $\omega(G)\ge 3$.
\end{theorem}

This bound is sharp for $\omega(G)=3$.  In 2007, Choudum, Karthick and Shalu
conjectured that $P_5$-free graphs have a quadratic binding function.

\begin{conj}{\em \cite{CKS07}}\label{Choudum06}
There is a constant $c$ such that $\chi(G)\le c\omega^2(G)$ if $G$ is $P_5$-free.
\end{conj}

Conjecture \ref{Choudum06} has been verified for many classes of $P_5$-free graphs, and tight linear binding functions are obtained for some  $(P_5, H)$-free graphs with $|V(H)|\le 5$, see \cite{BRSV13, CHM19, CHM21, ACTK20, CS19, CKMM20, CKS07, DXX2021, DX21+, FGMT95, SHTK19, S16}.
Very recently, Scott, Seymour and Spirkl \cite{SSS22} provided a near polynomial binding function for  $P_5$-free graphs stating that  $\chi(G)\leq \omega(G)^{log_{2}\omega(G)}$ if $G$ is $P_5$-free.

Let $F$ and $H$ be two graphs. We say that $F$ is a {\em blow up} of $H$ if $F$ can be obtained from $H$ by replacing each vertex with an independent set and then replacing each edge with a complete bipartite graph. A 5-{\em ring} is a blow up of a 5-hole. In \cite{sumner1} (see also \cite{ELMM13}), Sumner characterized the structure of $(P_5, K_3)$-free graphs.

\renewcommand{\baselinestretch}{1}
\begin{theorem}\label{sumner} {\em (\cite{sumner1})}
A connected $(P_5, K_3)$-free graph is either bipartite or a 5-ring.
\end{theorem}\renewcommand{\baselinestretch}{1.2}

By Theorems~\ref{thm-esperet} and \ref{sumner}, we have that each $(P_5, K_4)$-free graph is 5-colorable. The graph $K_1+(K_1\cup K_3)$ can be obtained from $K_4$ by adding a new vertex joining to one vertex of the $K_4$. So, $K_4$-free graphs must be $(K_1+(K_1\cup K_3))$-free. Motivated by Theorem~\ref{thm-esperet}, we study the chromatic number of $(K_1+(K_1\cup K_3))$-free graphs.
Among other results on the chromatic number of $P_5$-free graphs, we proved in  \cite{DXX2021} that if $G$ is  $(P_5, K_1+(K_1\cup K_3))$-free then $\chi(G)\le 3\omega(G)+11$.  In this paper, we present a characterization to $(P_5, K_1\cup K_3)$-free graphs, and prove that each $(P_5, K_1\cup K_3)$-free graph is $(2\omega(G)-1)$-colorable. Based on this, we get a tight upper bound for the chromatic number of $(P_5, K_1+(K_1\cup K_3))$-free graphs.

Before introducing the main results of this paper, we need some new notations. Let $v\in V(G)$, and let $X$ be a subset of $V(G)$. We use $N_X(v)$ to denote the set of neighbors of $v$ in $X$. We say that $v$ is \textit{complete} to $X$ if $N_X(v)=X$,  and say that $v$ is \textit{anticomplete} to $X$ if $N_X(v)=\emptyset$. For two subsets $X$ and $Y$ of $V(G)$, we say that $X$ is {\em complete} to $Y$ if each vertex of $X$ is complete to $Y$, say that $X$ is {\em anticomplete} to $Y$ if each vertex of $X$ is anticomplete to $Y$.

Let $d(v, X)=\min_{x\in X}d(v, x)$  and call $d(v, X)$  {\it the distance} of a vertex $v$ to a subset $X$. Let $i$ be a positive integer  and $N^i_G(X)=\{y\in V(G)\backslash X | d(y, X)=i\}$. We call $N^i_G(X)$ the \textit{$i$-neighborhood} of $X$ and simply write $N^1_G(X)$ as $N_G(X)$. If no confusion may occur, we write $N^i(X)$ instead of $N^i_G(X)$,  and $N^i(\{v\})$ is denoted by $N^i(v)$ for short. A set $D$ is said to be a {\em dominating set} of $G$ if $V(G)=D\cup N(D)$.

Suppose that  $C=v_1v_2v_3v_4v_5v_1$ is a 5-hole of $G$. Let $M(C)=V(G)\setminus (V(C)\cup N(C))$. For a subset $T\subseteq \{1,2,3,4,5\}$, we define
$$N_T(C) = \{x\;|\; x\in N(C), \mbox{ and } v_ix\in E(G) \mbox{ if and only if } i\in T\}.$$
It is easy to check that for $k\in \{1,2,3,4,5\}$ and $l=k+2$, $N_{\{k, k+2\}}(C)=N_{\{l, l+3\}}(C)$ and $N_{\{k, k+2, k+3\}}(C)=N_{\{l, l+1, l+3\}}(C)$, where the summation of subindex is taken modulo $5$ (in this paper, the summations of subindex are always taken modulo  some integer $h$  and we always set $h+1\equiv 1$).
We define
$${\cal N}^{(2)}(C)=\cup_{1\leq i\leq 5} N_{\{i, i+2\}}(C),$$
$${\cal N}^{(3)}(C)=\cup_{1\leq i\leq 5} (N_{\{i, i+1, i+2\}}(C)\cup N_{\{i, i+1, i+3\}}(C)),$$ and
$${\cal N}^{(4)}(C)=\cup_{1\leq i\leq 5} N_{\{i, i+1, i+2, i+3\}}(C).$$

Let  $C_1=x_1x_2x_3x_4x_5x_1$ and $C_2=y_1y_2y_3y_4y_5y_1$ be two disjoint 5-cycles. Let ${\cal F}$ be the graph obtained from $C_1\cup C_2$  by adding edges $\cup_{1\le i\le 5}\{x_iy_i, x_{i}y_{i+1}, x_iy_{i+3}\}$. It is easy to verify that each independent set of ${\cal F}$ has size at most 3, and  $\chi({\cal F})=4$.

\begin{figure}[htbp]
	\begin{center}
		\includegraphics[scale=0.5]{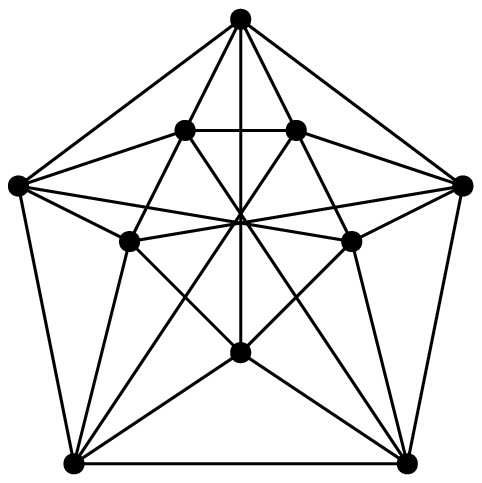}
	\end{center}
	\vskip -25pt
	{\small \caption{The graph ${\cal F}$}}
	\label{fig-1}
\end{figure}

The purpose of this paper is to prove the following

\renewcommand{\baselinestretch}{1}
\begin{theorem}\label{main2}
Let $G$ be  a connected $(P_5, K_1\cup K_3)$-free graph. Suppose that $G$ has non-dominating $5$-holes. Then, for each non-dominating $5$-hole $C=v_1v_2v_3v_4v_5v_1$, $V(G)$ can be partitioned into $4$ subsets $V(C)\cup {\cal N}^{(2)}(C)$, ${\cal N}^{(3)}(C)$, $N_{\{1,2,3,4,5\}}(C)$ and $M(C)$ with the following properties:
\begin{itemize}
\item [$(a)$] $G[V(C)\cup {\cal N}^{(2)}(C)]$ is a blow up of $C$, and $G[V(C)\cup {\cal N}^{(2)}(C)\cup {\cal N}^{(3)}(C)]$ is a blow up of a subgraph of ${\cal F}$,
\item [$(b)$] $M(C)\cup V(C)\cup {\cal N}^{(2)}(C)$ is complete to $N_{\{1,2,3,4,5\}}(C)$, and
\item [$(c)$] $M(C)$ is anticomplete to $V(C)\cup {\cal N}^{(2)}(C)$  but complete to ${\cal N}^{(3)}(C)$, and $M(C)$ is independent if ${\cal N}^{(3)}(C)\neq\emptyset$.

%\item [$(c)$] $A_4$ is complete to $A_1$ but anticomplete to $A_2$, and
\end{itemize}
\end{theorem}\renewcommand{\baselinestretch}{1.2}

\renewcommand{\baselinestretch}{1}
\begin{theorem}\label{main0}
If $G$ is  $(P_5, K_1\cup K_3)$-free then $\chi(G)\le 2\omega(G)-1$.
\end{theorem}

\begin{theorem}\label{main1}
If $G$ is a $(P_5, K_1+(K_1\cup K_3))$-free graph then  $\chi(G)\le \max\{2\omega(G), 15\}$, and there exists a $(P_5, K_1+(K_1\cup K_3))$-free graph $G$ with  $\chi(G)=2\omega(G)$.
\end{theorem}

\renewcommand{\baselinestretch}{1.2}

The proof of Theorem~\ref{main1} is heavily relied on  Theorem~\ref{main0}. The upper bound of Theorem~\ref{main0} is clearly tight  as $C_5$ and its blow up are extremal graphs.  We can construct a $(P_5, K_1+(K_1\cup K_3))$-free graph $G$ with  $\chi(G)=2\omega(G)$. Let $C=v_1v_2v_3v_4v_5v_1$ be a 5-hole. Let $H$ be the graph obtained from $C$ by replacing each vertex  $v_i$ of  $C$ by a 5-hole $C^i$,  for $1\le i\le 5$,  such that a vertex of $C^i$ and a vertex of $C^j$ are adjacent in $H$ if and only if $v_i$ is adjacent to $v_j$ in $C$.

It is certain that $H$ is $(P_5, K_1+(K_1\cup K_3))$-free and $\omega(H)=4$. We claim that $\chi(H)=8$. Without loss of generality, for each coloring  $\phi$ of $H$, we can always suppose  that $\phi(V(C^1))=\{1, 2, 3\}$ and $\phi(V(C^2))=\{4, 5, 6\}$. Let $\phi(V(C^3))=\{1, 7, 8\}$, $\phi(V(C^4))=\{3, 4, 5\}$ and $\phi(V(C^5))=\{6, 7, 8\}$. We see that  $\chi(H)\le 8$. If $\chi(H)\le 7$, then we may assume by symmetry that $\phi(V(C^3))=\{1,2,7\}$, but now we only have five colors $\{3, 4, 5, 6, 7\}$ that can be used to color $V(C^4)\cup V(C^5)$, a contradiction. Therefore, $\chi(H)=8=2\omega(H)$.

The following lemma, which is devoted to the structure of $P_5$-free graphs, will be used frequently in our proof. Here the summation of subindexes is taken modulo $5$.

\renewcommand{\baselinestretch}{1}
\begin{lemma}\label{P5free} {\em (\cite{DXX2021, ELMM13})}
Let $G$ be a $P_5$-free graph with a $5$-hole $C=v_1v_2v_3v_4v_5v_1$. Then
\begin{itemize}

\item [$(a)$] for $i\in \{1,2,3,4,5\}$, $N_{\{i\}}(C)=N_{\{i, i+1\}}(C)=\emptyset$, and $N_{\{i, i+2\}}(C)\cup N_{\{i, i+1, i+2\}}(C)$ is anticomplete to  $N^2(C)$,

\item [$(b)$] if $x\in N(C)$ and $N^2(x)\cap N^3(C)\neq\emptyset$ then $x\in N_{\{1,2,3,4,5\}}(C)$, and

\item [$(c)$] for each vertex $x\in N^2(C)$ and each component $B$ of $G[N^3(C)]$, $x$ is  either complete or anticomplete to  $B$.

\end{itemize}

\end{lemma}\renewcommand{\baselinestretch}{1.2}

The next section is devoted to the proofs of Theorem~\ref{main2} and Theorem~\ref{main0}. Theorem~\ref{main1} is proved in Sections 3.

\section{$(P_5, K_1\cup K_3)$-free graphs}

This section is aimed to prove Theorem~\ref{main2} and Theorem~\ref{main0}.
In this section, we always suppose that $G$ is a  $(P_5, K_1\cup K_3)$-free graph. If $G$ has a 5-hole, we always use $C=v_1v_2v_3v_4v_5v_1$ to denote a 5-hole in $G$. Recall that we define $M(C)=V(G)\setminus(V(C)\cup N(C))$.

\renewcommand{\baselinestretch}{1}
\begin{lemma}\label{lem-p5-K1-K3-structure}
If $G$ has a $5$-hole, then the followings hold for each $i\in \{1, 2, \ldots, 5\}$.
\begin{itemize}
	\item[$(a)$] $N_{\{i,i+1,i+2\}}(C)=N_{\{i,i+1,i+2,i+3\}}(C)=\emptyset$.

	\item[$(b)$] Both $N_{\{i,i+2\}}(C)$ and $N_{\{i, i+1, i+3\}}(C)$ are independent, and $N_{\{i,i+2\}}(C)$ is complete to  $N_{\{i+1, i+3\}}(C)\cup N_{\{i+1,i+4\}}(C)$.

\item[$(c)$] ${\cal N}^{(2)}(C)$ is complete to $N_{\{1,2,3,4,5\}}(C)$, and ${\cal N}^{(3)}(C)$ is complete to $M(C)$.

		\item[$(d)$] $N_{\{i, i+1, i+3\}}(C)$ is anticomplete to  $N_{\{i, i+3\}}(C)\cup N_{\{i+1, i+3\}}(C)$. Moreover,  if $M(C)\neq\emptyset$, then $N_{\{i, i+1, i+3\}}(C)$ is anticomplete to $N_{\{i, i+2, i+3\}}(C)\cup N_{\{i+1, i+3, i+4\}}(C)$, and is either complete or anticomplete to $N_{\{i-1, i, i+2\}}(C)\cup N_{\{i+1, i+2, i+4\}}(C)$ whenever both $N_{\{i-1, i, i+2\}}(C)$ and $N_{\{i+1, i+2, i+4\}}(C)$ are not empty.

		\item[$(e)$] If $\omega(G[N_{\{1,2,3,4,5\}}(C)])=\omega(G)-2$ or $M(C)\neq \emptyset$ then $G[V(C)\cup {\cal N}^{(2)}(C)]$  is a 5-ring.
	
\end{itemize}
\end{lemma}\renewcommand{\baselinestretch}{1.2}
\pf Suppose that  $N_{\{i,i+1,i+2\}}(C)\cup N_{\{i,i+1,i+2,i+3\}}(C)\neq \emptyset$ for some $i\in \{1, 2, \ldots, 5\}$. Let $v\in N_{\{i,i+1,i+2\}}(C)\cup N_{\{i,i+1,i+2,i+3\}}(C)$. Then $G[\{v, v_{i+1}, v_{i+2}, v_{i+4}\}]=K_1\cup K_3$. Hence $(a)$ holds.

Next we prove $(b)$. Suppose, for some $i$, $N_{\{i, i+2\}}(C)$ is not independent. Let $uv$ be an edge in $G[N_{\{i,i+2\}}(C)]$. Then  $G[\{u, v, v_{i}, v_{i+3}\}]=K_1\cup K_3$, a contradiction. Similarly, if $N_{\{i,i+1,i+3\}}(C)$ is not independent, let $uv$ be an edge of $G[N_{\{i,i+1,i+3\}}(C)]$, then  $G[\{u,v,v_{i},v_{i+2}\}]=K_1\cup K_3$, which leads to a contradiction. If $N_{\{i, i+2\}}(C)$ is not complete to  $N_{\{i+1, i+3\}}(C)$  for some $i$, choose $u\in N_{\{i,i+2\}}(C)$ and $v\in N_{\{i+1,i+3\}}(C)$ with $uv\not\in E(G)$, then $uv_iv_{i+4}v_{i+3}v$ is an induced $P_5$ of $G$. A similar contradiction occurs if $N_{\{i, i+2\}}(C)$ is not complete to  $N_{\{i+1, i+4\}}(C)$. Therefore, $(b)$ holds.

If  there exist $u\in N_{\{1,2,3,4,5\}}(C)$ and $v\in N_{\{i,i+2\}}$ with $uv\not\in E(G)$ for some $i$, then $G[\{u, v, v_{i+3}, v_{i+4}\}]=K_1\cup K_3$. If for some $i$,  there exist $u\in M(C)$ and $v\in N_{\{i, i+1,i+3\}}$ such that $uv\not\in E(G)$, then $G[\{u, v, v_{i}, v_{i+1}\}]=K_1\cup K_3$. This proves $(c)$.

If the first statement of $(d)$ is not true, then we may choose $u\in N_{\{i,i+1,i+3\}}(C)$ and $v\in N_{\{i, i+3\}}(C)\cup N_{\{i+1, i+3\}}(C)$ with $uv\in E(G)$ such that $G[\{u, v, v_{i}, v_{i+2}\}]=K_1\cup K_3$ when $v\in N_{\{i,i+3\}}(C)$, and  that $G[\{u, v, v_{i+1}, v_{i+4}\}]=K_1\cup K_3$ when $v\in N_{\{i+1, i+3\}}(C)$. Suppose that $M(C)\neq \emptyset$, and let $x\in M(C)$. Note that $M(C)$ is complete to ${\cal N}^{(3)}(C)$ by the statement $(c)$. If there exist $u\in N_{\{i,i+1,i+3\}}(C)$ and $v\in N_{\{i, i+2, i+3\}}(C)\cup N_{\{i+1, i+3, i+4\}}(C)$  with $uv\in E(G)$, then  $G[\{u, v, v_{i+4}, x\}]=K_1\cup K_3$ when $v\in N_{\{i, i+2, i+3\}}(C)$, and $G[\{u, v, v_{i+2}, x\}]=K_1\cup K_3$ when $v\in N_{\{i+1, i+3, i+4\}}(C)$.
Suppose that  $N_{\{i-1, i, i+2\}}(C)\neq\emptyset$ and $N_{\{i+1, i+2, i+4\}}(C)\neq\emptyset$. By symmetry, assume that $u\in N_{\{i, i+1, i+3\}}(C)$ is adjacent to $v\in N_{\{i-1, i, i+2\}}(C)$ and not adjacent to $w\in N_{\{i+1, i+2, i+4\}}(C)$. Then $vw\not\in E(G)$ and $G[\{u, v, v_i, w\}]=K_1\cup K_3$, a contradiction. Therefore, $(d)$ holds.

By the statement $(b)$, to prove that ${\cal N}^{(2)}(C)\cup V(C)$  induces a 5-ring, we only need to check that $N_{\{i, i+2\}}(C)$ is anticomplete to $N_{\{i, i+3\}}(C)\cup N_{\{i+2, i+4\}}(C)$.

By Lemma~\ref{P5free}$(a)$, we observe that $M(C)$ is anticomplete to ${\cal N}^{(2)}(C)\cup V(C)$. If $M(C)\ne \emptyset$,  then $N_{\{i, i+2\}}(C)$ must be anticomplete to $N_{\{i, i+3\}}(C)\cup N_{\{i+2, i+4\}}(C)$, otherwise a $K_1\cup K_3$ occurs.

Finally, suppose that $\omega(G[N_{\{1,2,3,4,5\}}(C)])=\omega(G)-2$, and let $K\subseteq  N_{\{1,2,3,4,5\}}(C)$ be a clique of size $\omega(G)-2$.  Assume by symmetry that $N_{\{i, i+2\}}(C)$ is not anticomplete to $N_{\{i+2, i+4\}}(C)$. Let $u\in N_{\{i,i+2\}}(C)$ and $v\in N_{\{i+2,i+4\}}(C)$ be an adjacent pair. Then  $K$ is complete to $\{u, v, v_{i+2}\}$ by statement $(c)$, and so $G$ contains a  clique of size $\omega(G)+1$. This leads to a contradiction and completes the proof of Lemma~\ref{lem-p5-K1-K3-structure}. \qed

\medskip

From Lemma~\ref{lem-p5-K1-K3-structure}$(a)$, we observe that
$$N(C)=N_{\{1,2,3,4,5\}}(C)\cup {\cal N}^{(2)}(C)\cup {\cal N}^{(3)}(C),$$
and it follows from Lemma~\ref{lem-p5-K1-K3-structure}$(d)$ that if $M(C)\ne \emptyset$ and  $N_{\{i,i+1,i+3\}}(C)\ne \emptyset$ for each $1\le i\le 5$,  then ${\cal N}^{(3)}(C)$ is either independent or induces a 5-ring in $G$.

\medskip

\noindent{\bf Proof of Theorem}~\ref{main2}: Suppose that $G$ has a non-dominating 5-hole $C=v_1v_2v_3v_4v_5v_1$, that is, $M(C)=V(G)\setminus(V(C)\cup N(C))\neq\emptyset$. Let $A_1=V(C)\cup {\cal N}^{(2)}(C)$, $A_2={\cal N}^{(3)}(C)$, and $A_3=N_{\{1,2,3,4,5\}}(C)$.

By Lemma~\ref{lem-p5-K1-K3-structure}$(b)$ and $(d)$,  we observe that $G[A_1]$ is a 5-ring which is a blow up of $C$, and $G[A_1\cup A_2]$ is a blow up of a subgraph of ${\cal F}$.

By Lemma~\ref{lem-p5-K1-K3-structure}$(c)$, we have that $A_1$ is complete to $A_3$. To prove the second statement, we only need to verify that $A_3$ is complete to $M(C)$. If it is not the case, choose $u\in A_3$ and $v\in M(C)$ with $uv\not\in E(G)$,  then $G[\{u, v, v_1, v_2\}]=K_1\cup K_3$, a contradiction. Therefore, $(b)$ is true.

By Lemma~\ref{lem-p5-K1-K3-structure}$(c)$, we observe that $M(C)$ is anticomplete to $A_1$ and complete to $A_2$. Suppose that $A_2\neq \emptyset$, and let $x\in N_{\{i, i+1, i+3\}}(C)$ for some $i\in \{1, 2, 3, 4, 5\}$. If the third statement is not true then there must be an edge $uv$ in $G[M(C)]$ and so $G[\{u, v, v_{i+2}, x\}]=K_1\cup K_3$. This leads to a contradiction and proves $(c)$, and also completes the proof of Theorem~\ref{main2}. \qed

\bigskip

Now we turn to prove Theorem~\ref{main0}. The following two colorings will be used in the proof of Theorem~\ref{main0}.

By Lemma~\ref{lem-p5-K1-K3-structure}$(d)$, we can construct a 5-coloring $\psi$ of $G[V(C)\cup N(C)]-N_{\{1,2,3,4,5\}}(C)$ as below:
\begin{equation}\label{eqa-array-1}
\left\{\begin{array}{ll}
\psi^{-1}(1) = N_{\{2,4\}}(C)\cup  N_{\{1,2,4\}}(C)\cup \{v_3\}, \\
\psi^{-1}(2) =  N_{\{3,5\}}(C)\cup  N_{\{2,3,5\}}(C) \{v_4\},\\
\psi^{-1}(3) =  N_{\{1,4\}}(C)\cup N_{\{3,4,1\}}(C)\cup \{v_2, v_5\},\\
\psi^{-1}(4) = N_{\{2,5\}}(C)\cup N_{\{4,5,2\}}(C)\cup\{v_1\},\\
\psi^{-1}(5) =  N_{\{1,3\}}(C)\cup  N_{\{5,1,3\}}(C).
\end{array}\right.
\end{equation}

If  $M(C)\ne \emptyset$, it follows from Theorem~\ref{main2} that we can construct a 4-coloring $\phi$ of $G[V(C)\cup N(C)] - N_{\{1,2,3,4,5\}}(C)$ as below:
\begin{equation}\label{eqa-array-0}
\left\{\begin{array}{ll}
\phi^{-1}(1) =  N_{\{1,4\}}(C)\cup  N_{\{1, 2, 4\}}(C)\cup N_{\{2, 4\}}(C)\cup \{v_3,v_5\},\\
\phi^{-1}(2) =  N_{\{2, 5\}}(C)\cup  N_{\{2, 3, 5\}}(C)\cup N_{\{3,5\}}(C)\cup \{v_1,v_4\},\\
\phi^{-1}(3) = N_{\{1,3\}}(C)\cup  N_{\{1,3,4\}}(C)\cup N_{\{5, 1, 3\}}(C)\cup \{v_2\}, \\
\phi^{-1}(4) = N_{\{4, 5, 2\}}(C).
\end{array}\right.
\end{equation}

\noindent\textbf{Proof of Theorem~\ref{main0}.} Let $G$ be a connected $(P_5, K_1\cup K_3)$-free graph with $\omega(G)=h$. Clearly the theorem holds when $h=1$. If $h=2$, then $G$ is bipartite or a 5-ring by Theorem~\ref{sumner} and so $\chi(G)\le 3=2h-1$. Thus we assume that $h\geq 3$ and the theorem holds for all graphs with clique number smaller than $h$.

If $G$ does not have any 5-hole, then for an arbitrary vertex $v$, $G-N(v)$ is bipartite and $\omega(G[N(v)])\leq h-1$. Thus by induction $\chi(G)\leq 2+\chi(G[N(v)])\leq 2+(2(h-1)-1)=2h-1$. Otherwise, let $C=v_1v_2v_3v_4v_5v_1$ be a 5-hole of $G$.
It is certain that  $\omega(G[N_{\{1,2,3,4,5\}}(C)\cup N_{\{i, i+1, i+3\}}(C)])\le h-2$ for each $i\in \{1, 2, 3, 4, 5\}$ as $\{v_i, v_{i+1}\}$ is complete to $N_{\{1,2,3,4,5\}}(C)\cup N_{\{i, i+1, i+3\}}(C)$.

If $N_{\{1,2,3,4,5\}}(C)=M(C)=\emptyset$, then $\chi(G)\le 5\le 2h-1$ by the coloring $\psi$ defined in (\ref{eqa-array-1}).

If $N_{\{1,2,3,4,5\}}(C)=\emptyset$ and $M(C)\ne \emptyset$, then ${\cal N}^{(3)}(C)\neq\emptyset$, which implies that $M(C)$ is independent by  Theorem~\ref{main2}$(c)$. It follows from the coloring $\phi$ defined in (\ref{eqa-array-0}) that $\chi(G)\le 5\le 2h-1$.

Suppose that $N_{1,2,3,4,5}(C)\neq \emptyset$ and $M(C)= \emptyset$.
If $\omega(G[N_{\{1,2,3,4,5\}}(C)])\le h-3$, then $\chi(G-N_{\{1,2,3,4,5\}}(C))\le 5$ by the coloring $\psi$ defined in (\ref{eqa-array-1}), which implies that $\chi(G)\le \chi(G-N_{\{1,2,3,4,5\}}(C))+\chi(G[N_{\{1,2,3,4,5\}}(C)])\le 5+(2(h-3)-1)<2h-1$ by induction. So, suppose that $\omega(G[N_{\{1,2,3,4,5\}}(C)])= h-2$. By Lemma~\ref{lem-p5-K1-K3-structure}$(d)$ and $(e)$, we have that $N_{\{1, 3\}}(C)$ is anticomplete to $N_{\{1,4\}}(C)\cup N_{\{3,4,1\}}(C)\cup \{v_2, v_5\}$, and so we can modify the coloring $\psi$ by recoloring  $N_{\{1, 3\}}(C)$ with 3, which implies that $\chi(G-N_{\{1,2,3,4,5\}}(C)\cup N_{\{5, 1, 3\}}(C))\le 4$. Now, we have that $\chi(G)\le \chi(G-N_{\{1,2,3,4,5\}}(C)\cup N_{\{5, 1, 3\}}(C))+\chi(G[N_{\{1,2,3,4,5\}}(C)\cup N_{\{5, 1, 3\}}(C)])\le 4+(2(h-2)-1)=2h-1$ by induction.

Therefore, suppose that $N_{\{1,2,3,4,5\}}(C)\ne \emptyset$ and $M(C)\ne \emptyset$. Thus,  ${\cal N}^{(2)}(C)\cup V(C)$  induces a 5-ring by Lemma~\ref{lem-p5-K1-K3-structure}$(d)$. It is obvious that $G[M(C)]$ is $K_3$-free, otherwise a triangle of $G[M(C)]$ together with any vertex of $C$ induces a $K_1\cup K_3$, and so $\chi(G[M(C)])\le 3$ by Theorem~\ref{sumner}.

If ${\cal N}^{(3)}(C)=\emptyset$, then $\chi(G-N_{\{1,2,3,4,5\}}(C))=3$ as $M(C)$ is anticomplete to ${\cal N}^{(2)}(C)$ by Lemma~\ref{P5free}$(a)$, and so $\chi(G)\le 3+(2(h-2)-1)=2h-1$ by induction.
Thus, suppose that ${\cal N}^{(3)}(C)\ne\emptyset$, which implies that $M(C)$ is independent by Theorem~\ref{main2}$(c)$. By the coloring $\phi$ defined in (\ref{eqa-array-0}),
$G-N_{\{1,2,3,4,5\}}(C)\cup N_{\{4, 5, 1\}}(C)\cup M(C)$ is 3-colorable, and so $G-N_{\{1,2,3,4,5\}}(C)\cup N_{\{4, 5, 1\}}(C)$ is 4-colorable.
Since $\omega(G[N_{\{1,2,3,4,5\}}(C)\cup N_{\{4, 5, 1\}}(C)])\le h-2$, we have $\chi(G)\le 4 + \chi(G[N_{\{1,2,3,4,5\}}(C)\cup N_{\{4, 5, 1\}}(C)])\le 4+(2(h-2)-1)=2h-1$ by induction.  This completes the proof of Theorem~\ref{main0}. \qed

\section{$(P_5, K_1+(K_1\cup K_3))$-free graphs}

Before proving Theorem~\ref{main1}, we first present several lemmas on the structure of $(P_5, K_1+(K_1\cup K_3))$-free graphs. From now on, we always suppose that $G$ is a connected $(P_5, K_1+(K_1\cup K_3))$-free graph without clique cutset. For two subsets $X$ and $Y$ of $V(G)$, we say that  $X$ is {\em adjacent} to $Y$ if $N(X)\cap Y\neq\emptyset$.

Let $C=v_1v_2v_3v_4v_5v_1$ be a 5-hole of $G$. Recall that  ${\cal N}^{(2)}(C)=\cup_{1\le i\le 5} N_{\{i, i+2\}}(C)$, ${\cal N}^{(4)}(C)=\cup_{1\le i\le 5} N_{\{i,i+1, i+2, i+3\}}(C)$, and $M(C)=V(G)\setminus(V(C)\cup N(C))$. We further define ${\cal N}^{(3, 1)}(C)=\cup_{1\le i\le 5} N_{\{i,i+1, i+2\}}(C)$, and ${\cal N}^{(3, 2)}(C)=\cup_{1\le i\le 5} N_{\{i,i+1, i+3\}}(C)$.
By Lemma~\ref{P5free}$(a)$, we have
$$N(C)=N_{\{1,2,3,4,5\}}(C)\cup {\cal N}^{(2)}(C)\cup {\cal N}^{(3, 1)}(C)\cup {\cal N}^{(3, 2)}(C)\cup {\cal N}^{(4)}(C).$$

\renewcommand{\baselinestretch}{1}
\begin{lemma}{\em (\cite{DXX2021})}\label{P5K1K3K1free}
Let $C=v_1v_2v_3v_4v_5v_1$ be a  $5$-hole of $G$,  and $T$ be a component of $G[N^2(C)]$. Then the followings hold.
\begin{itemize}

\item [$(a)$] For each $i\in \{1,2,3,4,5\}$, $G[N(v_i)]$ is $(K_1\cup K_3)$-free, $G[N_{\{i, i+2\}}(C)]$ is $K_3$-free,  and $N_{\{i, i+1, i+2\}}(C)\cup N_{\{i,i+1,i+3\}}(C)\cup N_{\{i,i+1,i+2,i+3\}}(C)$ is independent.

\item [$(b)$] If no vertex in $N(C)$ dominates $T$, then there exist two non-adjacent vertices $u$ and $v$ in $N(C)$ such that both $N_T(u)$ and $N_T(v)$ are not empty.

\end{itemize}
\end{lemma}

\renewcommand{\baselinestretch}{1}
\begin{lemma}\label{P5K1K3K1freeC5-1}
Let  $C=v_1v_2v_3v_4v_5v_1$ be a $5$-hole of $G$, $S$ be a component of $G[N_{\{1,2,3,4,5\}}(C)]$ with $\omega(S)\ge 2$. Then for each $i\in \{1,2,3,4,5\}$, the followings hold.
\begin{itemize}

\item [$(a)$] $N_{\{i,i+2\}}(C)\cup N_{\{i,i+1,i+2\}}(C) $ is complete to $S$, and $N_{\{i, i+2\}}(C)$ is independent.

\item [$(b)$] For each edge $xy$ in $S$, no vertex of $N_{\{i,i+1,i+3\}}(C)\cup N_{\{i,i+1,i+2,i+3\}}(C)$ is anticomplete to $\{x, y\}$.

\item [$(c)$] $N_{\{i, i+2\}}(C)$ is anticomplete to $N_{\{i-1,i,i+1\}}(C)\cup N_{\{i-1,i,i+2\}}(C)\cup N_{\{i-1,i,i+1,i+2\}}(C)$.

\item[$(d)$] $\chi(G - N_{\{1,2,3,4,5\}}(C)-M(C))\leq 5$.

\end{itemize}
\end{lemma}\renewcommand{\baselinestretch}{1.2}
\pf Suppose that, for some $i\in \{1,2,3,4,5\}$, $N_{\{i,i+2\}}(C)\cup N_{\{i,i+1,i+2\}}(C)$ has a vertex $u$ that is not complete to  $S$. If $u$ is anticomplete to $S$, then $G[\{u,v,w,v_i,v_{i+4}\}]=K_1+(K_1\cup K_3)$. Otherwise, there exists an edge, say $vw$ in $S$, such that $uv\in E(G)$ and $uw\notin E(G)$. Then  $G[\{u,v,w, v_{i+3},v_{i+4}\}]=K_1+(K_1\cup K_3)$. Both are contradictions.

Suppose that $N_{\{i,i+2\}}(C)$ is not independent for some $i\in \{1,2,3,4,5\}$. Choose an edge  $xy$ in  $G[N_{\{i,i+2\}}(C)]$, and let $z$ be an arbitrary vertex  of $S$. Then $xz\in E(G)$ and $yz\in E(G)$, and so $G[\{x,y,z, v_{i},v_{i+3}\}]=K_1+(K_1\cup K_3)$, a contradiction. Therefore, $(a)$ holds.

Let $xy$ be an edge of $S$,   and $v\in N_{\{i,i+1,i+3\}}(C)\cup N_{\{i,i+1,i+2,i+3\}}(C)$. If $vx\notin E(G)$ and $vy\notin E(G)$, then $G[\{v, x, y, v_{i+3}, v_{i+4}\}]=K_1+(K_1\cup K_3)$.
Therefore, $(b)$ holds.

Suppose that $(c)$ is  not true for some $i\in \{1,2,3,4,5\}$. Let $v\in  N_{\{i, i+2\}}(C)$ and $u\in  N_{\{i-1,i,i+1\}}(C)\cup N_{\{i-1,i,i+2\}}(C)\cup  N_{\{i-1,i,i+1,i+2\}}(C)$ such that $uv\in E(G)$. By $(a)$ and $(b)$, we observe that there exists a vertex $w\in N_{\{1,2,3,4,5\}}(C)$ such that $wu\in E(G)$ and $wv\in E(G)$, which implies that $G[\{v,u,w, v_{i},v_{i+3}\}]=K_1+(K_1\cup K_3)$. Therefore, $(c)$ is true.

By $(a)$, $(c)$ and Lemma~\ref{P5K1K3K1free}$(a)$, we have that $N_{\{i,i+2\}}(C)\cup N_{\{i-1, i, i+1\}}(C)\cup N_{\{i-1,i,i+2\}}(C)\cup N_{\{i-1,i,i+1,i+2\}}(C)$ is independent for each $i\in \{1,2,3,4,5\}$. By coloring $\{v_{i+3}\}\cup N_{\{i,i+2\}}(C)\cup N_{\{i-1, i, i+1\}}(C)\cup N_{\{i-1,i,i+2\}}(C)\cup N_{\{i-1,i,i+1,i+2\}}(C)$ with color $i$, we get a $5$-coloring of
$G - N_{\{1,2,3,4,5\}}(C)-M(C)$. This proves $(d)$, and completes the proof of Lemma~\ref{P5K1K3K1freeC5-1}. \qed

\renewcommand{\baselinestretch}{1}
\begin{lemma}\cite{DXX2021}\label{P5K1K3K1freeC5}
Let $C=v_1v_2v_3v_4v_5v_1$ be a $5$-hole of $G$. Then $G[N^3(C)]$ is $K_3$-free, and  $N^2(C)$ can be partition into two parts $A$ and $B$  such that both  $G[A]$ and $G[B]$ are $K_3$-free.
\end{lemma}

\begin{lemma} \label{P5K1K3K1freeC5-2}
Let $C=v_1v_2v_3v_4v_5v_1$ be a $5$-hole of $G$, and $S$ be a component of $G[N_{\{1,2,3,4,5\}}(C)]$. If $N(S)\cap N^2(C)\neq \emptyset$, then $N(x)\cap N^2(C)=N(y)\cap N^2(C)$ for any $x, y\in S$.
\end{lemma}\renewcommand{\baselinestretch}{1.2}
\pf Suppose that $N(S)\cap N^2(C)\neq \emptyset$. We apply induction on $|S|$.  The lemma holds trivially if $|S|=1$. Suppose that $|S|=k\ge 2$, and the lemma holds on all components of $G[N_{\{1,2,3,4,5\}}(C)]$ of size less than $k$. There must be a vertex, say $x$, in $S$ such that $S-x$ is connected, and $N(S-x)\cap N^2(C)\neq \emptyset$.  Let $y$ be a neighbor of $x$ in $S$. To prove the lemma, we only need to verify that $N(x)\cap N^2(C)=N(y)\cap N^2(C)$. Suppose that it is not the case. Then, we may assume, without loss of generality, that $u\in N(x)\cap N^2(C)$  and $u\notin N(y)\cap N^2(C)$, which implies that $G[\{x,y,u,v_1,v_{2}\}]=K_1+(K_1\cup K_3)$. This leads to a contradiction and proves the lemma. \qed

\renewcommand{\baselinestretch}{1}
\begin{lemma} \label{P5K1K3K1freeC5-4}
Let $C=v_1v_2v_3v_4v_5v_1$ be a $5$-hole of $G$, and $T$ be a component of $G[N^2(C)]$. Suppose that ${\cal N}^{(3,2)}(C)\cup {\cal N}^{(4)}(C)\neq\emptyset$. Then
\begin{itemize}
\item [$(a)$] $T$ is a single vertex adjacent to $N_{\{1,2,3,4,5\}}(C)$ if $N(T)\cap ({\cal N}^{(3,2)}(C)\cup {\cal N}^{(4)}(C))\ne \emptyset$  and $\omega(G[N_{\{1,2,3,4,5\}}(C)])\geq 2$, and
\item [$(b)$] $T$ is $K_3$-free if $N(T)\cap ({\cal N}^{(3,2)}(C)\cup {\cal N}^{(4)}(C))=\emptyset$.
\end{itemize}
\end{lemma}\renewcommand{\baselinestretch}{1.2}
\pf Let $Q={\cal N}^{(3,2)}(C)\cup {\cal N}^{(4)}(C)$.

Firstly, we prove $(a)$. Suppose that $N(T)\cap Q\ne\emptyset$ and $\omega(G[N_{\{1,2,3,4,5\}}(C)])\geq 2$.  Since $N(T)\cap Q\ne\emptyset$, we have that, for some $i\in \{1,2,3,4,5\}$, $N_{\{i, i+1, i+3\}}\cup N_{\{i, i+1, i+2, i+3\}}$ has a vertex $u$ that is complete to $T$, otherwise an induced $P_5$ appears in $G$. By Lemma~\ref{P5K1K3K1freeC5-1}$(b)$, $u$ has a neighbor, say $v$, in $N_{\{1,2,3,4,5\}}(C)$. If $v$ is anticomplete to $T$ then $G[\{u, v, v_i, v_{i+1}, z\}=K_1+(K_1\cup K_3)$ for any vertex $z\in T$. This proves that  $N(T)\cap N_{\{1,2,3,4,5\}}(C)\ne\emptyset$, that is, $T$ is adjacent to $N_{\{1,2,3,4,5\}}(C)$.

Suppose that $|V(T)|\ge 2$, and let $xy$ be an edge of $T$. Since $G$ is $P_5$-free, we have that, for some $i\in \{1,2,3,4,5\}$,  $N_{\{i,i+1,,i+3\}}(C)\cup N_{\{i,i+1,i+2,i+3\}}(C)$  has a vertex, say $u$, that is complete to $T$. Particularly, $\{ux, uy\}\subseteq E(G)$. By Lemma \ref{P5K1K3K1freeC5-1}$(b)$, we may choose a neighbor $v$ of $u$ in $N_{\{1,2,3,4,5\}}(C)$. If $\{vx, vy\}\subseteq E(G)$ then $G[\{u,v,v_{i+4},x,y\}=K_1+(K_1\cup K_3)$. Otherwise, we may assume by symmetry that $vx\not\in E(G)$, then $G[\{u, v, v_i, v_{i+1}, x\}]=K_1+(K_1\cup K_3)$. Therefore, $(a)$ holds.

Suppose to the contrary of $(b)$ that $N(T)\cap Q=\emptyset$ and  $T$ has a $K_3$, say $w_1w_2w_3w_1$. Let $u$ be a vertex in $Q$, and suppose that $uv_1, uv_2\in E(G)$ by symmetry. Since ${\cal N}^{(3,1)}(C)$ is anticomplete to $N^2(C)$ by Lemma~\ref{P5free}$(a)$, we may choose a vertex, say $x$, in $N_{\{1,2,3,4,5\}}(C)$, and let $x'$ be a neighbor of $x$ in $T$. To avoid a $K_1+(K_1\cup K_3)$ on $\{u, v_1, v_2, x, x'\}$, we have that $ux\not\in E(G)$. If $x$ is complete to $T$, then $G[\{v_1, w_1, w_2, w_3, x\}]=K_1+(K_1\cup K_3)$. Otherwise, there must be an edge $y_1y_2$ in $T$ such that $xy_1\in E(G)$ and $xy_2\not\in E(G)$, and so $uv_1xy_1y_2$ is an induced $P_5$. This proves $(b)$ and  Lemma~\ref{P5K1K3K1freeC5-4}. \qed

\medskip

Let $A$ be an  antihole  with $V(A)=\{v_1, v_2, \cdots, v_k\}$. We enumerate the vertices of $A$ cyclically  such that $v_iv_{i+1}\notin E(G)$  and simply write $A=v_1v_2\cdots v_k$. Here the summations of subindex are taken modulo $k$  and we set $k+1\equiv 1$.

Suppose that $G$ induces an antihole $A=v_1v_2\cdots v_{k}$ with $k\geq 6$. We use $S(A)$ to denote the set of vertices which are complete to $A$, and let $T(A)=N(A)\backslash S(A)$. Note that $T(A)$ is not complete to $A$. For each $i\in \{1, 2, \ldots, k\}$, we define $T_i(A)$ to be the subset of $T(A)$ such that for each vertex $x$ of $T_i(A)$, $i$ is the minimum index with $xv_i\in E(G)$ and $xv_{i-1}\not\in E(G)$.

Clearly, $T(A)=\cup_{1\le i\le k} T_i(A)$, and $T_i(A)\cap T_j(A)=\emptyset$ if $i\ne j$. Since $G$ is $K_1+(K_1\cup K_3)$-free, we have that $G[S(A)]$ is $K_1\cup K_3$-free, and   $G[T_i(A)]$ is $K_1\cup K_3$-free for each $i\in \{1, 2, \ldots, k\}$.

The following lemma was proved in \cite{DXX2021} without using the notations  $T_i(A)$. Here we present its short proof.

\renewcommand{\baselinestretch}{1}
\begin{lemma}\label{P5K1K3K1freeantihole}
Let $G$ be a $(P_5, C_5, K_1+(K_1\cup K_3))$-free graph, $A=v_1v_2\cdots v_{k}$ an antihole of $G$ with $k\ge 6$. Then $T_i(A)$ is independent for each $i\in\{1,2,\ldots, k\}$, and $N^2(A)=\emptyset$.
\end{lemma}\renewcommand{\baselinestretch}{1.2}
\pf Let $i\in \{1,2, \ldots, k\}$. Firstly,  for each vertex $v$ of $T_i(A)$,
\begin{equation}\label{eqa-v-vi+2-2}
vv_{i+2}\in E(G),
\end{equation}
as otherwise  either $vv_iv_{i+2}v_{i-1}v_{i+1}$ is an induced $P_5$ when $vv_{i+1}\not\in E(G)$ or $vv_iv_{i+2}v_{i-1}v_{i+1}v$ is a 5-hole when $vv_{i+1}\in E(G)$.

Suppose that $T_i(A)$ is not independent. Let $x$ and $x'$ be two adjacent vertices of $T_i(A)$. Then  $G[\{v_{i-1}, v_{i}, v_{i+2}, x, x'\}]=K_1+(K_1\cup K_3)$ by (\ref{eqa-v-vi+2-2}). Therefore, $T_i(A)$ is an independent set.

Suppose that $N^2(A)\neq \emptyset$. Let $v$ be a vertex in $N(A)$ that has a neighbor, say $x$, in $N^2(A)$. If $v\in S(A)$ then $G[\{v, v_1, v_3, v_5, x\}]= K_1+(K_1\cup K_3)$. Otherwise, we may assume that $v\in T_1(A)$ by symmetry. By (\ref{eqa-v-vi+2-2}), $G[\{v, v_1, v_{3}, v_{5}, x\}]=K_1+(K_1\cup K_3)$ if $vv_{5}\in E(G)$, and a $P_5=xvv_1v_{5}v_{2k+1}$ appears if $vv_{5}\not\in E(G)$. Therefore, $N^2(A)=\emptyset$. \qed

\medskip

\noindent\textbf{Proof of Theorem~\ref{main1}.}  Let $G$ be a $\{P_5, K_1+(K_1\cup K_3\}$-free graph with $\omega(G)=h$. We may suppose that $G$ is connected, contains no clique cutset, and is not perfect. Thus, $h\ge 2$ as  $G$ must induce a 5-hole  or an odd antihole with at least 7 vertices by the Strong Perfect Graph Theorem \cite{CRST06}.

When $h\in\{2, 3\}$, the theorem follows immediately from Theorems~\ref{thm-esperet} and \ref{sumner}. Suppose that $h\ge 4$, and the theorem holds for all $\{P_5, K_1+(K_1\cup K_3\}$-free graphs with clique number less than $h$.

Since $G$ is $P_5$-free, it is certain that $N^4(S)=\emptyset$ for any subset $S$ of $V(G)$.

Let $\gamma=2h-5$. We distinguish two cases depending on the existence of $5$-holes in $G$, and will use two color sets ${\cal C}_1=\{\alpha_1, \alpha_2, \alpha_3, \alpha_4, \alpha_5\}$ and ${\cal C}_2=\{\beta_1, \beta_2,\cdots, \beta_{\gamma}\}$ to color $G$.

Firstly, suppose that $G$ induces no 5-holes. Then, $G$ must induce an antihole of size at least 6. Let $A=v_1v_2\cdots v_{k}$, where $k\ge 6$, be an antihole of $G$. Let $S$ be the set of vertices that are complete to $A$, and let $T=V(G)\backslash(A\cup S)$. It is clear that $G[S]$ is $K_1\cup K_3$-free. By Lemma~\ref{P5K1K3K1freeantihole}, $V(G)=A\cup S\cup T$.

For integer $i\in \{1, 2, \ldots, k\}$, let $T_i$ be the subset of $T$ such that for each vertex $x$ of $T_i$, $i$ is the minimum index with $xv_i\in E(G)$ and $xv_{i-1}\not\in E(G)$.
By Lemma~\ref{P5K1K3K1freeantihole}, $T_i\cup\{v_{i-1}\}$ is an independent set.

If $S\neq \emptyset$, then $\chi(G[A \cup T])\leq  k$ by Lemma~\ref{P5K1K3K1freeantihole}, and so $\chi(G) \leq k + (2(h-\lfloor{k\over 2}\rfloor)-1) \le 2h$ by induction. Therefore, we suppose that $S=\emptyset$.

We further suppose that $A$ has the least number of vertices under the assumption that $k\geq 6$.   Notices that $\frac{k}{2}\leq h$ if $k$ is even and  $\frac{k-1}{2}\leq h$ if $k$ is odd.
If $k\le 15$, then $\chi(G)\leq k\le 15$ by Lemma~\ref{P5K1K3K1freeantihole}.
If $h>\lfloor\frac{k}{2}\rfloor$ then $\chi(G)\leq 2\lceil\frac{k}{2}\rceil\leq 2h$. So, we suppose that $h=\lfloor\frac{k}{2}\rfloor\ge 8$.

Since $S=\emptyset$, for each vertex $v\in T$, there must exist an integer $i$ such that  $vv_i\in E(G)$ and $vv_{i-1}\not\in E(G)$. For integer $i\in \{1, 2, \ldots, k\}$, let $T_i$ be the subset of $T$ such that for each vertex $x$ of $T_i$, $i$ is the minimum index with $xv_i\in E(G)$ and $xv_{i-1}\not\in E(G)$. By Lemma~\ref{P5K1K3K1freeantihole}, $T_i\cup\{v_{i-1}\}$ is an independent set.

If $k$ is even, then by coloring the vertices in $T_i\cup\{v_{i-1}\}$ with color $i$, we get a $2h$-coloring of $G$. Therefore, we suppose that $k$ is odd.

Let $v$ be a vertex in $T_i$ for some $i$.

If $vv_{i+2}\not\in E(G)$, then $G[\{v, v_{i-1}, v_{i}, v_{i+1}, v_{i+2}\}]$ is a $C_5$ or $P_5$ depending on $vv_{i+1}\in E(G)$ or not. So, $vv_{i+2}\in E(G)$.   We will show that
\begin{equation}\label{eqa-vi-vi-1}
\mbox{if $\{vv_i, vv_{i+2}\}\subseteq E(G)$, then $vv_{i+1}\in E(G)$ and $vv_{i-2}\not\in E(G)$.}
\end{equation}
First suppose $vv_{i+1}\notin E(G)$. If $vv_{i+4}\in E(G)$, then $G[\{v, v_i, v_{i+1}, v_{i+2}, v_{i+4}\}]=K_1+(K_1\cup K_3)$.  If $vv_{i+4}\notin E(G)$, then $G[\{v, v_{i+1}, v_{i+2}, v_{i+3}, v_{i+4}\}]$ is a $C_5$ or $P_5$ depending on $vv_{i+3}\in E(G)$ or not. Both are contradictions. This shows that $vv_{i+1}\in E(G)$.  If $vv_{i-2}\in E(G)$, then $G[\{v, v_{i-2}, v_{i-1}, v_i, v_{i+2}\}]=K_1+(K_1\cup K_3)$. Therefore, (\ref{eqa-vi-vi-1}) holds.

Without loss of generality, we suppose that $\{vv_1, vv_2, vv_3\}\subseteq E(G)$ and $vv_{k}\notin E(G)$. By (\ref{eqa-vi-vi-1}) and by symmetry, $vv_{k-1}\notin E(G)$.

If $vv_4\notin E(G)$, then  $vv_5\notin E(G)$ by (\ref{eqa-vi-vi-1}) and by symmetry, and hence  $vv_j\in E(G)$ for all $j\in \{6, \ldots, k-2\}$ to avoid an induced $K_1+(K_1\cup K_3)$ on $\{v, v_2, v_4, v_j, v_k\}$.  But then  $G[\{v, v_5, v_6, \ldots, v_{k-1}\}]$ is an antihole with less vertices, which contradicts the choice of $A$. Therefore, $vv_4\in E(G)$.

If $vv_5\notin E(G)$, then $G[\{v, v_1, v_2, v_3, v_4, v_5, v_{k}\}]$ is an antihole with less vertices, which is contradiction to the choice of $A$. So, $vv_5\in E(G)$. Repeating this argument, we have $vv_j\in E(G)$ for all $j\in \{1, 2, \ldots, k-3\}$.

If $vv_{k-2}\in E(G)$, then we have a clique of size at leat $\lfloor\frac{k}{2}\rfloor+1$, which contradicts $h=\lfloor\frac{k}{2}\rfloor$. Thus, $vv_{k-2}\notin E(G)$. Consequently, we have, by symmetry, that each vertex in $T$ is nonadjacent to exactly three consecutive vertices of $A$.

Suppose that there exist $x_1\in T_1$ and $x_k\in T_k$ with $x_1x_k\in E(G)$. Then, $\{x_1, x_k\}$ is complete to $V(A)\setminus\{v_k, v_{k-1}, v_{k-2}, v_{k-3}\}$. Since $k$ is odd, we have that $\{v_1, v_3, \ldots, v_{k-4}\}$ induces a $K_{k-3\over 2}$, which together with $\{x_1, x_k\}$ induces a $K_{h+1}$. Therefore, we may suppose by symmetry that $T_i$ is anticomplete to $T_{i+1}$ for any $i\in \{1, 2, \ldots, k\}$. Thus, $G$ is a subgraph of a blow up of $A$ by Lemma~\ref{P5K1K3K1freeantihole}, which implies $\chi(G)=h+1<2h$.

This shows that Theorem~\ref{main1} holds if $G$ does not induce 5-holes. From now on to the end of this paper, we always assume that $G$ induces a 5-hole, and
\begin{equation}\label{eqa-minimum-omega}
\mbox{let $C=v_1v_2v_3v_4v_5v_1$ be a 5-hole of $G$ that minimizes $\omega(G[N_{\{1,2,3,4,5\}}(C)])$.}
\end{equation}

Recall that we can partition $N(C)$ into 5 subsets: ${\cal N}^{(2)}=\cup_{1\le i\le 5}N_{\{i,i+2\}}(C)$, ${\cal N}^{(3,1)}=\cup_{1\le i\le 5}N_{\{i,i+1,i+2\}}(C)$,
${\cal N}^{(3,2)}=\cup_{1\le i\le 5} N_{\{i,i+1,i+3\}}(C)$,
${\cal N}^{(4)}=\cup_{1\le i\le 5} N_{\{i,i+1,i+2,i+3\}}(C)$,
and $N_{\{1,2,3,4,5\}}(C)$.

By Lemma~\ref{P5K1K3K1freeC5}, $N^3(C)$ is $K_3$-free, and $N^2(C)$ can be partitioned into two subsets each of which induces a $K_3$-free subgraph.
Thus by Theorem~\ref{sumner}, we have that $\chi(G[N^2(C)])\le 6$ and $\chi(G[N^3(C)])\le 3$.

By Lemma~\ref{P5K1K3K1free}$(a)$, we have that, for each $i\in \{1, 2, 3, 4, 5\}$, $G[N_{\{i, i+2\}}(C)]$ is $K_3$-free,  and $\{v_{i+4}\}\cup N_{\{i, i+1, i+2\}}(C)\cup N_{\{i,i+1,i+3\}}(C)\cup N_{\{i,i+1,i+2,i+3\}}(C)$ is independent. If $G[N_{1, 3}(C)\cup N_{1, 4}(C)]$ is not $K_3$-free, let $xyzx$ be a triangle in $G[N_{1, 3}(C)\cup N_{1, 4}(C)]$, then $G[\{v_1, v_5, x, y, z\}]=K_1+(K_1\cup K_3)$, a contradiction. So, we have by symmetry that $G[N_{1, 3}(C)\cup N_{1, 4}(C)]$ and $G[N_{2, 4}(C)\cup N_{2, 5}(C)]$ are both $K_3$-free. Hence we may conclude that $\chi(G[V(C)\cup {\cal N}^{(3,1)}\cup {\cal N}^{(3,2)}\cup {\cal N}^{(4)} \cup N^3(C)])\leq 5$, and $\chi(G[{\cal N}^{(2)}\cup N^2(C)])\leq 9$
as ${\cal N}^{(2)}$ is anticomplete to $N^2(C)$ by Lemma~\ref{P5free}.

If $N_{\{1,2,3,4,5\}}(C)$ is independent, then $\chi(G)\leq \chi(G[V(C)\cup {\cal N}^{(3,1)}\cup {\cal N}^{(3,2)}\cup {\cal N}^{(4)} \cup N^3(C)])+\chi(G[{\cal N}^{(2)}\cup N^2(C)])+\chi(G[N_{\{1,2,3,4,5\}}(C)])\le 5+9+1=15$.  Hence, we may assume that
$$2\le \omega(G[N_{\{1,2,3,4,5\}}(C)])\le \omega(G)-2=h-2.$$

Let $Q={\cal N}^{(3, 1)}\cup {\cal N}^{(3, 2)}\cup {\cal N}^{(4)}$.

We claim that if $Q\neq\emptyset$ then
\begin{equation}\label{eqa-Q-N2C-1}
\mbox{$N^2(C)$ is anticomplete to each non-isolated component of $G[N_{\{1,2,3,4,5\}}(C)]$.}
\end{equation}
If it is not the case, then let $xy$ be an edge of some non-isolated component of $G[N_{\{1,2,3,4,5\}}(C)]$. By Lemma~\ref{P5K1K3K1freeC5-2}, $N^2(C)$ has a vertex, say $u$, complete to $\{x, y\}$. By Lemma~\ref{P5K1K3K1freeC5-1}$(b)$ and by symmetry, $Q$ has a vertex, say $v$, adjacent to $x$. Without loss of generality, we may assume that $\{vv_1, vv_2\}\subseteq E(G)$. Thus $G[\{u, v, v_1, v_2, x\}]=K_1+(K_1\cup K_3)$, a contradiction. Therefore, (\ref{eqa-Q-N2C-1}) holds.

\subsection{Suppose that $2\le \omega(G[N_{\{1,2,3,4,5\}}(C)])\leq h-3$}

In this case, we have that $h\ge 5$. Let $\omega(G[N_{\{1,2,3,4,5\}}(C)])=t$. Note that  by Lemma~\ref{P5K1K3K1freeC5-1} $\chi(G-N_{\{1,2,3,4,5\}}(C)-N^{2}(C)-N^{3}(C))\le 5$, and by Theorem~\ref{main0} $\chi(G[N_{\{1,2,3,4,5\}}(C)])\le 2t-1$  as $G[N_{\{1,2,3,4,5\}}(C)]$ is $K_1\cup K_3$-free.

If $N^2(C)=\emptyset$, then $\chi(G)\leq 5+(2t-1)< 2h$ by induction. Thus we may assume that $N^2(C)\neq\emptyset$, and  without loss of generality, $G[N^2(C)]$ is connected.

Recall that $\chi(G[N_{\{1,2,3,4,5\}}(C)])\leq 2h-7$ by induction, and $\chi(G[N^2(C)\cup V(C)\cup {\cal N}^{(2)}])\leq 6$ by Lemma~\ref{P5K1K3K1freeC5}.

If $Q=\emptyset$, then color $V(C)\cup {\cal N}^{(2)}\cup N^2(C)$ with ${\cal C}_1\cup\{\beta_1\}$ and color $N_{\{1,2,3,4,5\}}(C)\cup N^3(C)$ with ${\cal C}_2\setminus\{\beta_1\}$. Thus, we obtain a $2h$-coloring of $G$.

Therefore, we further suppose that $Q\neq\emptyset$.

Let $N^{2,0}(C)\subseteq N^2(C)$ be the set of vertices anticomplete to $Q$. If $Q$ is anticomplete to $N^2(C)$, that is, $N^2(C)=N^{2,0}(C)$,  then $N^2(C)$ is anticomplete to all non-isolated components of $G[N_{\{1,2,3,4,5\}}(C)]$ by (\ref{eqa-Q-N2C-1}), which implies that $N^3(C)=\emptyset$. We can color $V(C)\cup N(C)$ with ${\cal C}_1\cup {\cal C}_2$ such that all isolated vertices of $G[N_{\{1,2,3,4,5\}}(C)]$ receive the same color $\beta_1$, and color $N^2(C)$ with ${\cal C}_1\cup {\cal C}_2\setminus\{\beta_1\}$ (this is certainly reasonable as $\chi(G[N^2(C)])\leq 6$ by Lemma~\ref{P5K1K3K1freeC5}).

Suppose that $Q$ is adjacent to $N^2(C)$. By Lemma~\ref{P5K1K3K1freeC5-4}, each vertex of $N^2(C)\setminus N^{2,0}(C)$ is an isolated component of $G[N^2(C)]$. Since $N^{2, 0}(C)$ is anticomplete to $Q\cup (N^2(C)\setminus N^{2, 0}(C))$, by Lemma~\ref{P5K1K3K1freeC5-1}$(d)$ and Lemma~\ref{P5K1K3K1freeC5}, we can color $G-N_{\{1,2,3,4,5\}}(C)-N^3(C)$ with ${\cal C}_1\cup \{\beta_1\}$. Since $\omega(G[N_{\{1,2,3,4,5\}}(C)])\le h-3$ and $G[N^3(C)]$ is $K_3$-free, we can color $N_{\{1,2,3,4,5\}}(C)\cup N^3(C)$ with ${\cal C}_2\setminus \{\beta_1\}$ by Theorems~\ref{sumner} and \ref{main0}. Therefore, $\chi(G)\le \chi(G-N_{\{1,2,3,4,5\}}(C)-N^3(C))+ \chi(G[N_{\{1,2,3,4,5\}}(C)\cup N^3(C)])\le 2h$. Thus when $2\le \omega(G[N_{\{1,2,3,4,5\}}(C)])\leq h-3$, $\chi(G)\le 2h$.

\subsection{Suppose that $\omega(G[N_{\{1,2,3,4,5\}}(C)])= h-2$}

Now, suppose that $\omega(G[N_{\{1,2,3,4,5\}}(C)])=h-2$.  By (\ref{eqa-minimum-omega}), we have that
\begin{equation}\label{eqa-omega-1}
\mbox{$\omega(G[N_{\{1,2,3,4,5\}}(C')])=h-2$ for each 5-hole $C'$ of $G$.}
\end{equation}

Let $S$ be a component of $G[N_{\{1,2,3,4,5\}}(C)]$ with $\omega(S)=h-2$.

By Lemma~\ref{P5K1K3K1freeC5-1}$(a)$,
${\cal N}^{(2)}\cup {\cal N}^{(3,1)}$ is complete to $S$. Hence we have that
\begin{equation}\label{m1}
\mbox{${\cal N}^{(3,1)}=\emptyset$ and $V(C)\cup {\cal N}^{(2)}(C)$ induces a 5-ring}
\end{equation}
as otherwise we can find a clique of size at least $\omega(G)+1$.

By Lemma~\ref{P5K1K3K1freeC5-1}$(d)$, we can define a 5-coloring $\phi$ on $G-N_{\{1,2,3,4,5\}}(C)-N^2(C)-N^3(C)$  with color set ${\cal C}_1$ as followsing:

\begin{equation}\label{eqa-3-array-0}
\left\{\begin{array}{ll}
\phi^{-1}(\alpha_1)=\{v_1\}\cup N_{\{3,5\}}(C)\cup N_{\{2,3,5\}}(C)\cup N_{\{2,3,4,5\}}(C) \\
\phi^{-1}(\alpha_2)=\{v_2\}(C)\cup N_{\{4,1\}}(C)\cup N_{\{3,4,1\}}(C)\cup N_{\{3,4,5,1\}}(C)\\
\phi^{-1}(\alpha_3)=\{v_3\}\cup N_{\{5,2\}}(C)\cup N_{\{4,5,2\}}(C)\cup N_{\{4,5,1,2\}}(C) \\
\phi^{-1}(\alpha_4)=\{v_4\}\cup N_{\{1,3\}}(C)\cup N_{\{5,1,3\}}(C)\cup N_{\{5,1,2,3\}}(C)\\ \phi^{-1}(\alpha_5)=\{v_5\}\cup N_{\{2,4\}}(C)\cup N_{\{1,2,4\}}(C)\cup N_{\{1,2,3,4\}}(C).
\end{array}\right.
\end{equation}

If $N^2(C)=\emptyset$, then by Theorem~\ref{main0}, $\chi(G)\leq 5+2(h-2)-1=2h$  as $G[N_{\{1,2,3,4,5\}}(C)]$ is $K_1\cup K_3$-free.

Thus suppose that $N^2(C)\ne \emptyset$, and without loss of generality, suppose that $G[N^2(C)]$ is connected.

By (\ref{m1}), we have that ${\cal N}^{(3,1)}=\emptyset$, and so $Q={\cal N}^{(3,2)}\cup {\cal N}^{(4)}$. Let $N^{2,0}(C)\subseteq N^2(C)$ be the set of vertices anticomplete to $Q$.

We first suppose that $Q\ne\emptyset$, and discuss two cases depending upon whether  $N^2(C)$ is adjacent to $Q$.

\noindent{\bf Case} 1. Suppose that $Q$ is anticomplete to $N^2(C)$. Then each component of $G[N^{2}(C)]$  is $K_3$-free by Lemma~\ref{P5K1K3K1freeC5-4}, and $N^2(C)$ is anticomplete to all non-isolated components of $G[N_{\{1,2,3,4,5\}}(C)]$ by (\ref{eqa-Q-N2C-1}).
Consequently we have that  $G[N_{\{1,2,3,4,5\}}(C)]$ has isolated components (as $N^2(C)\neq \emptyset$) and also has non-isolated components (as $\omega(G[N_{\{1,2,3,4,5\}}(C)])=h-2\ge 2$). If $N^3(C)\neq \emptyset$, let $n_3\in N^3(C)$, $n_2\in N^2(C)$ be a neighbor of $n_3$, $s_1$ an isolated component of $G[N_{\{1,2,3,4,5\}}(C)]$ with $s_1n_2\in E(G)$, and $s_2s'_2$ be an edge of some component of $G[N_{\{1,2,3,4,5\}}(C)]$, then $n_3n_2s_1v_1s_2$ is an induced $P_5$, a contradiction. Therefore, $N^3(C)=\emptyset$.

Now, we can color  $G[N_{\{1,2,3,4,5\}}(C)]$ with color set ${\cal C}_2$ such that such that all isolated vertices of $G[N_{\{1,2,3,4,5\}}(C)]$ receive the same color $\beta_1\in {\cal C}_2$, and color $G[N^2(C)]$ with the colors in ${\cal C}_1\cup {\cal C}_2\setminus\{\beta_1\}$ (this is reasonable as $\chi(G[N^2(C)])\leq 6$ by Lemma~\ref{P5K1K3K1freeC5}). This together with the 5-coloring defined in (\ref{eqa-3-array-0}) gives a $2h$-coloring of $G$.

\medskip

\noindent{\bf Case} 2. Suppose that $N^2(C)$ is adjacent to $Q$. By Lemma~\ref{P5K1K3K1freeC5-4}, we have that each component of $G[N^{2, 0}(C)]$ is $K_3$-free, and each of the other components of $G[N^2(C)]$ is a single vertex. Since  $\omega(G[N_{\{1,2,3,4,5\}}(C)\cup N_{\{5,1,3\}}(C)\cup N_{\{5,1,2,3\}}(C)])=h-2$, we have that $\chi(G[N_{\{1,2,3,4,5\}}(C)\cup N_{\{5,1,3\}}(C)\cup N_{\{5,1,2,3\}}(C)])\le 2h-5$ by induction. Using the 5-coloring $\phi$ defined in (\ref{eqa-3-array-0}), we can construct a 5-coloring of  $G[V(C)\cup N^2(C)\cup (N(C)\setminus (N_{\{1,2,3,4,5\}}(C)\cup N_{\{5,1,3\}}(C)\cup N_{\{5,1,2,3\}}(C))]$ by coloring all the vertices of $N^2(C)\setminus N^{2,0}(C)$ by $\alpha_4$, and coloring all the vertices of $N^{2,0}(C)$ by $\{\alpha_1, \alpha_2, \alpha_3\}$ (this is reasonable by Lemma~\ref{P5K1K3K1freeC5-4}). Then by  coloring $N^3(C)$ with 3 colors used on $G[N_{\{1,2,3,4,5\}}(C)\cup N_{\{5,1,3\}}(C)\cup N_{\{5,1,2,3\}}(C)]$, we have that $\chi(G)\leq 5+(2h-5)=2h$ by induction.

\medskip

We have shown that $\chi(G)\leq 2h$ when $Q\ne\emptyset$. Next, we suppose that $Q=\emptyset$.

If $N^2(C)$ is adjacent to only isolated component of  $G[N_{\{1,2,3,4,5\}}(C)]$, we see that $N^3(C)=\emptyset$ by the same argument as that used in Case 1, then we can color $G[N_{\{1,2,3,4,5\}}(C)]$ with color set ${\cal C}_2$  such that all isolated components receive $\beta_1$, and color $N^2(C)$ with ${\cal C}_1\cup {\cal C}_2\setminus\{\beta_1\}$ (this reasonable as $\chi(G[N^2(C)])\leq 6$ by Lemma~\ref{P5K1K3K1freeC5}. This together with $\phi$ defined in (\ref{eqa-3-array-0}) is certainly a $2h$-coloring of $G$.

So, we suppose that $N^2(C)$ is  adjacent to some non-isolated  components of $G[N_{\{1,2,3,4,5\}}(C)]$, and let $S_1$ be the vertex set of such a component. Let  $S_2=N_{\{1,2,3,4,5\}}(C)\backslash S_1$, $T_1=N(S_1)\cap N^2(C)$,  and $T_2=N^2(C)\backslash T_1$. It is obvious that $S_1$ is anticomplete to $T_2$, and is complete to $T_1$ by Lemma~\ref{P5K1K3K1freeC5-2}.

Therefore, $G[T_1]$ is $K_3$-free. Note that $G[N^2(C)]$ is connected by our assumption. To avoid an induced $P_5$ starting from $T_2$ and terminating on $C$, each component of $G[T_2]$ is dominated by some vertex of $T_1$, and consequently $G[T_2]$ is $K_3$-free too. We will show that \begin{equation}\label{eqa-T2-independent}
\mbox{$T_2$ is independent.}
\end{equation}

If it is not the case, let $Z$ be a non-isolated component of $G[T_2]$, let $t_1\in T_1$ be a vertex complete to $Z$, and $s_2\in S_2$ be a  vertex adjacent to $Z$. If $s_2$ is not complete to $Z$, let $z_1z_2$ be an edge of $Z$ such that $s_2z_1\in E(G)$ and $s_2z_2\not\in E(G)$, then $z_2z_1s_2v_1s_1$ is an induced $P_5$ for any vertex $s_1\in S_1$, a contradiction. Therefore, $s_2$ is complete to $Z$. If $s_2t_1\not\in E(G)$, then for any vertices $s_1\in S_1$ and $z\in V(Z)$, $C'=s_1t_1zs_2v_1s_1$ is a 5-hole with $N_{\{1,2,3,4,5\}}(C')=\emptyset$, a contradiction to (\ref{eqa-omega-1}). So, we have further  that $s_2t_1\in E(G)$. But now, we have a $K_1+(K_1\cup K_3)$ induced by $\{s_2, t_1, v_1\}$ together with any two adjacent vertices of $Z$. Therefore, (\ref{eqa-T2-independent}) holds.

Note that $G[N_{\{1,2,3,4,5\}}(C)]$ is $(K_1\cup K_3)$-free and $G[N^3(C)]$ is $K_3$-free by Lemmas~\ref{P5K1K3K1free} and \ref{P5K1K3K1freeC5}, we see that $\chi(G[N_{\{1,2,3,4,5\}}(C)\cup N^3(C)])\leq 2h-5$ by Theorems~\ref{sumner} and \ref{main0}. Since $T_2$ is independent by (\ref{eqa-T2-independent}), we have that $\chi(G[N^2(C)\cup {\cal N}^{(2)}(C)\cup V(C)])\leq5$,  and so $\chi(G)\leq 2h$ as desired.  This completes the proof of Subsection 3.2,  and also proves Theorem~\ref{main1}. \qed

\noindent{\bf Acknowledgement: We thank Dr. Karthick for pointing out an error in our earlier version on the construction of some extremal graphs}.

\section{Statements and Declarations}

\subsection {Funding}
This work was supported by National Natural Science Foundation of China (No. 11931106 and 12101117) and by Natural Science Foundation of Jiangsu Province (No. BK20200344). Author Baogang Xu has received research support from National Natural Science Foundation of China. Author Yian Xu has received research support from National Natural Science Foundation of China and Natural Science Foundation of Jiangsu Province. 

\subsection{Competing Interests}
The authors have no relevant financial or non-financial interests to disclose.

\subsection{Author Contributions}

All authors contributed to the study conception and design. Material preparation, data collection and analysis were performed by Wei Dong, Baogang Xu and Yian Xu. The first draft of the manuscript was written by Wei Dong, Baogang Xu and Yian Xu, and all authors commented on previous versions of the manuscript. All authors read and approved the final manuscript.

\section{Data Availability Statements}

No data applicable. 

\end{document}